\documentstyle[12pt,saks]{article}
\def\cc{\setcounter{equation}{0}}
\def\be{\begin{equation}}
\def\ee{\end{equation}}

\newcommand{\blem}{\begin{lemma}}
\newcommand{\elem}{\end{lemma}}
\newcommand{\bth}{\begin{theo}}
\newcommand{\eth}{\end{theo}}
\newcommand{\bcor}{\begin{corol}}
\newcommand{\ecor}{\end{corol}}
\newcommand{\br}{\begin{rem}}
\newcommand{\er}{\end{rem}}
\newcommand{\bprob}{\begin{prob}}
\newcommand{\eprob}{\end{prob}}
\newcommand{\ba}{\begin{array}}
\newcommand{\ea}{\end{array}}
\newcommand{\beq}{\begin{eqnarray}}
\newcommand{\beqq}{\begin{eqnarray*}}
\newcommand{\eeq}{\end{eqnarray}}
\newcommand{\eeqq}{\end{eqnarray*}}

\newcounter{minutes}
\divide\time by 60
\newcounter{hours}
\multiply\time by 60
\addtocounter{minutes}{-\time}

\begin{document}

\bigskip
\bigskip

\begin{center}
{\Large\bf On quasiplanes in Euclidean spaces }
\end{center}

\medskip

\begin{center}
{\large\bf O.~Martio, V.M.~Miklyukov, and M.~Vuorinen} 
\texttt{FILE:~\jobname .tex, 2005-11-06, 
        printed: \number\year-\number\month-\number\day, 
        \thehours.\ifnum\theminutes<10{0}\fi\theminutes}
\end{center}

\vspace{2mm}

{\em $2000$ Mathematics Subject Classification.}  30C62, 30C65


\bigskip

{\bf Abstract.} A variational inequality for the images of $k$-dimensional
hyperplanes
under quasiconformal maps of the $n$-dimensional Euclidean space
is proved when  $1\le k\le n-2 \, .$
\bigskip

\section{Main Results}{}

Below we use the terminology and notation of \cite{FMMVW}.

Let ${\bf R}^n$ be the $n-$dimensional Euclidean space, $n\ge 2$,
let $f:{\bf R}^n\rightarrow {\bf R}^n$ be a mapping
of the class $W^1_{n,\,\rm loc}({\bf R}^n)$, and
let $f':{\bf R}^n\to {\bf R}^n$ denote its formal derivative.
We write
$$
\|f'(x)\|=\max_{|h|=1}|f'(x)\,h|\,.
$$

A homeomorphism $f:{\bf R}^n\rightarrow {\bf R}^n$ is called
{\it $K$-quasiconformal} \cite[p.~250-252]{HKM} if
$f\in W^1_{n,\,\rm loc}({\bf R}^n)$ and
\be\label{eq1}
\|f'(x)\|^n\leq K\, J(x,f) \quad \mbox{a.e. on} \quad  {\bf R}^n,
\ee
where $J(x,f)={\rm det}\,(f'(x))$ is the Jacobian of $f$ at the point $x\in
{\bf R}^n.$
The smallest constant $K$ in (\ref{eq1}) is called the outer
dilatation of $f$ and denoted by $K_O(f)$. The smallest constant $K\ge 1$ in
the inequality
$$
J(x,f)\le K\, \ell (f'(x))^n
$$
is called the inner dilatation of $f:{\bf R}^n\to{\bf R}^n$ and denoted
by $K_I(f)$. Here, $ \ell(f'(x))=\min_{|h|=1}|f'(x)\,h|\,.$
The quantity
$$
K(f)=\max\{K_O(f),K_I(f)\}
$$
is called the maximal dilatation of $f$ \cite[Section {\bf 14.1}]{HKM}.

If $f:{\bf R}^n\rightarrow {\bf R}^n$ is quasiconformal, then it is
well-known that $f({\bf R}^n) = {\bf R}^n $ and the
inverse map $f^{-1}:{\bf R}^n\rightarrow {\bf R}^n$
is also quasiconformal in ${\bf R}^n$ with $K(f^{-1})=K(f) \,.$

Let $d(x',x'')$ be the Euclidean distance for $x', x'' \in {\bf R}^n \,.$
We write
$$
B(a,t)=\{x\in{\bf R}^n:d(a,x)<t\},\quad S(a,t)=\{x\in{\bf R}^n:d(a,x)=t\}\,.
$$

Let $y=(y_1,y_2,...,y_n)$ be a point in ${\bf R}^n$ and $1\le k\le n-2$.
Consider a $k$-dimensional plane
$$
\Pi^k_0=\{y=(y_1,y_2,\ldots ,y_n)\in {\bf R}^n:\ \ y_{k+1}=\ldots=y_n=0\}\,.
$$
Let $\Pi^k$ be a $k$-dimensional surface in ${\bf R}^n$, i.e. a
homeomorphic image of a $k$-dimensional plane. Fix a point $a\in \Pi^k$ and
$R,\,0<R<\infty$. Denote by $\Pi^k(a,R)$ a component of
$\Pi^k\cap B(a,R)$, $a\in \Pi^k(a,R)$. We say that a $k$-dimensional
surface $\Pi^k$ is a {\it $K$-quasiplane},
if there exists
a $K$-quasiconformal mapping $f:{\bf R}^n\to {\bf R}^n$ such that
$f(\Pi^k)= \Pi^k_0$.

For $n=2$, the following result is well known \cite{A}.
{\it A curve $\Pi^1\subset {\bf R}^2$ is a $K$-quasiconformal line
if and only if there exists a constant $C(K)$ such that
$$
 \frac{d(\zeta_2,\zeta_1)}{d(\zeta_3,\zeta_1)}\leq C(K)
$$
for every triple of distinct points on $\Pi^1$ such that $\zeta_2$ lies between
$\zeta_1$ and $\zeta_3$.}

For more information on quasiconformal lines, we refer to \cite{G},
\cite[Chapter 14]{AVV}, \cite{Gr}, \cite{Kr}, \cite{MMV3},
\cite{MMPV}, \cite{Vu2}. Quasiplanes $\Pi^{n-1}$ in ${\bf R}^n$
with ${\rm codim}\,\Pi^{n-1}=1$
were considered in \cite{MV}, \cite{VVW} and, in the setup of
Riemannian manifolds, in  \cite{Mik}.
Quasiplanes may have a highly complicated structure. For instance,
it is well-known that quasiconformal lines may have Hausdorff dimension
$>1$ and hence be non-rectifiable; on the other hand sufficient
conditions for the $(n-1)$-rectifiability of quasiplanes
 $\Pi^{n-1}$ in ${\bf R}^n$ were given in \cite{MV}.
Below we consider the case of $K$-quasiplanes in $ {\bf R}^n$ of dimension
$1\le k\le n-2$.

For an arbitrary open subset $\Sigma $ on the sphere
$ S(a,r)\subset {\bf R}^n$ we define the quantity
\begin{equation}
\label{f11}
\eta(\Sigma)=\sup_{A}\inf\limits_{\varphi}\frac{\left(\displaystyle
\int_{\Sigma}|\nabla_{S}\varphi|^n d{\cal
H}^{n-1}\right)^{1/n}}{\left(\displaystyle
\int\limits_{\Sigma}|\varphi-A|^n d {\cal
H}^{n-1}\right)^{1/n}}\,.
\end{equation}

Here the infimum is taken over all functions
$$
\varphi\in W^1_n(\Sigma),\quad \varphi|_{\partial\Sigma}=0\,,
$$
the supremum is taken over all constants $A$ and
$\partial \Sigma$ is the boundary with respect to $S(a,r) \,.$

The quantity $\nabla_S\varphi$ is the gradient of $\varphi$ on $S(a,r)$.
(The most convenient for our purposes the definition of $\nabla_S\varphi$ on
a surface $S\subset {\bf R}^n$ see, for example, 
in \S 2 of \cite{Simons:1968}).

The zero boundary values
$ \varphi|_{\partial\Sigma}=0\,$ are understood in the Sobolev sense, i.e.
$\varphi\in W^1_{n, 0}(\Sigma) \,.$

It is clear that
\begin{equation}
\label{eqezo}
\eta(\Sigma)\ge \inf\limits_{\varphi}\frac{\left(\displaystyle
\int_{\Sigma}|\nabla_{S}\varphi|^n
d{\cal H}^{n-1}\right)^{1/n}}{\left(\displaystyle
\int\limits_{\Sigma}|\varphi|^n
d {\cal H}^{n-1}\right)^{1/n}}\,.
\end{equation}
Because ${\rm dim}\,S=n-1\, ,$  it follows from Sobolev's
$C^0$-embedding theorem that the right side of (\ref{eqezo}) is $> 0$
at least for every $\Sigma$ with a smooth boundary
$\partial\Sigma\ne \emptyset \,.$
To find the best $\eta(\Sigma)$ is an open problem.
\medskip

\begin{theo}{}
\label{theoknl}
Let $\Pi^k$ be a $K$-quasiplane in ${\bf R}^n$ with $1\le k \le n-2$.
Then for every point $a\in \Pi^k$ and all
$0<r<R<\infty$, the following relation holds
\begin{equation}
\label{eqth1}
{\rm exp}\left\{{{1}\over{K_O(f)}}\int\limits_{r}^{R}
\eta(\Sigma(a,\tau))\,d\tau\right\}\,\le
\,D(n,K)\,\left(\frac{R}{r} \right)^{n\beta}\,.
\end{equation}
Here $\beta= K^{1/(n-1)}$ and
$$
D(n,K)=D^{2n}_*\,,\quad
D_*= D_*(K)=\, {\mathrm{exp}} \, (4K(K+1)\sqrt{K-1}) \, .
$$
\end{theo}
\bigskip

\cc
\section{The language of differential forms}{}

Let $1\le k \le n-2$ and let
$$
f=(f_1,\ldots,f_k,f_{k+1},\ldots,f_n):{\bf R}^n \to {\bf R}^n
$$
be a $K$-quasiconformal mapping.

Let $\Pi^k_0$ be a $k$-dimensional plane of the form
$$
\Pi^k_0=\{y=(y_1,\ldots,y_k,y_{k+1},\ldots,y_n):y_{k+1}=\ldots=y_n=0\}\,.
$$
We will study some properties of the $k$-dimensional surface
$$
\Pi^k=f^{-1}(\Pi^k_0)\,.
$$

Fix a point $a\in \Pi^k$. For $R >0$ let
$\Pi(a,R)$ be the component of $\Pi^k\cap B(a,R)$,
containing the point $a$ and
$\Sigma(a,R)=S(a,R)\setminus \overline{\Pi(a,R)}$.

We consider the differential form of degree $n-1$
$$
\omega=\sum_{i=k+1}^{n}(-1)^{i-k-1}f_i\,df_{1}\wedge\ldots\wedge \widehat{df_i}
\wedge\ldots\wedge df_n\,.
$$
It is clear that for a.e. $R\in (0,\infty)$, we have
$$
\omega\in W^{1,1}(S(a,R))\quad\mbox{and}\quad \omega|_{\partial
\Sigma(a,R)}=0\,.
$$
We have
$$
\begin{array}{ll}
d\omega&=\sum_{i=k+1}^n(-1)^{i-k-1}df_i\,\wedge df_{k+1}\wedge\ldots\wedge \widehat{df_i}\wedge
\ldots\wedge df_n\\ \\
\quad&=(n-k)df_1\wedge\ldots\wedge df_n\, ,\\ \\
\end{array}
$$
and thus,
$$
*d\omega=(n-k)\,J(x,f)\ge 0\,.
$$

For a.e. $R\in (0,\infty)$ we write
$$
\varepsilon_f(a,R)=\sup_{\omega_0}\int\limits_{\Sigma(a,R)}
(*d\omega)\,d{\cal H}^{n-1}
\left/\int\limits_{\Sigma(a,R)}|\omega-\omega_0|\,d{\cal H}^{n-1}\right.\,,
$$
where the supremum is taken over all weakly closed $(n-1)$-forms
$\omega_0\in W^{1,1}(B(a,2R))$ such that
$$
\omega_0 \left|_{S(a,R)}\in W^{1,1}(S(a,R)) \right.  \, .
$$
Denote
$$
V(a,r)=\int\limits_{B(a,r)}J(x,f)\,d{\cal H}^n\,.
$$
Because ${\cal H}^n(\Pi_0)=0$ we see that ${\cal H}^n(\Pi^k)=0$
and for a.e. $r>0$ we have
$$
{\cal H}^{n-1}(\Pi^k\cap S(a,r))=0\,.
$$
Thus,
$$
\int\limits_{\Sigma(a,r)}\omega_0=\int\limits_{S(a,r)\setminus \Pi^k}\omega_0=
\int\limits_{S(a,r)}\omega_0\,.
$$

For an arbitrary weakly closed \cite{FMMVW} form $\omega_0$ and almost all
$0<R<\infty$ the following relations hold
\begin{equation}
\begin{array}{ll}
\displaystyle\int\limits_{\Sigma(a,R)}|\omega-\omega_0|\,d{\cal H}^{n-1}&\ge
\displaystyle\int\limits_{\Sigma(a,R)}(\omega-\omega_0)=\\ \\
\quad&=(n-k)\displaystyle\int\limits_{B(a,R)}J(x,f)\,d{\cal H}^{n}=(n-k)V(a,R)\,.\\ \\
\end{array}
\label{eq1840}
\end{equation}
Here we used that $f_i|_{\Pi^k}=0$ $(i=k+1,\ldots,n)$ and 
that for the weakly closed form $\omega_0$:
$$
\begin{array}{ll}
\displaystyle\int\limits_{\Sigma(a,R)}(\omega-\omega_0)&=
\displaystyle\int\limits_{\Sigma(a,R)}\omega=\\ \\
\quad&=\displaystyle\int\limits_{\Sigma(a,R)}\sum_{i=k+1}^{n}(-1)^{i-k-1}f_i\,df_{1}\wedge\ldots\wedge \widehat{df_i}
\wedge\ldots\wedge df_n=\\ \\
\quad&=\sum_{i=k+1}^{n}(-1)^{i-k-1}
\displaystyle\int\limits_{\Sigma(a,R)}f_i\,df_{1}\wedge\ldots\wedge \widehat{df_i}
\wedge\ldots\wedge df_n=\\ \\
\quad&=\sum_{i=k+1}^{n}\displaystyle\int\limits_{B(a,R)}
df_{1}\wedge\ldots\wedge df_i \wedge\ldots\wedge df_n=\\ \\
\quad&=(n-k)\displaystyle\int\limits_{B(a,R)}J(x,f)\,d{\cal H}^{n}\,.\\ \\
\end{array}
$$
The relation
\begin{equation}
\label{eqOGF}
\int\limits_{\Sigma(a,R)}f_i\,df_{1}\wedge\ldots\wedge \widehat{df_i}
\wedge\ldots\wedge df_n=(-1)^{i-k-1}
\int\limits_{B(a,R)}
df_{1}\wedge\ldots\wedge df_n
\end{equation}
is clear in the case of $C^2$-smooth maps. In the general case we approximate
the homeomorphism $f:{\bf R}^n\to{\bf R}^n$ of class $W^{1,n}_{\rm loc}({\bf R}^n)$
by smooth maps with the following properties:

$i)\;$ $f_s\to f$ locally uniformly in ${\bf R}^n$ as $s\to\infty$;

$ii)\;$ $\|f'_s(x) - f'(x)\|_{L^n(D)}\to 0$ as $s\to\infty$ for every
subdomain $D\subset\subset {\bf R}^n$.

We easily see that such an approximation is possible if we use the familiar
technique of approximating Sobolev functions by smoothed averages (see, for
example, \cite[\S 4.2.1]{EG}). Then we fix arbitrary $r_1,\,r_2$ with
$0<r_1<R<r_2<\infty$. Using the formula (\ref{eqOGF})
for smooth functions $f_s$, we get
$$
(-1)^{i-k-1}\int\limits_{r_1}^{r_2}dr\int\limits_{B(a,r)}J(x,f_s)\,d{\cal H}^n=
\int\limits_{r_1}^{r_2}\int\limits_{\Sigma(a,R)}f_{si}\,df_{s1}\wedge\ldots\wedge
\widehat{df_{si}}\wedge\ldots\wedge df_{sn}\,.
$$
Since $f_s\to f$ uniformly on $B(a,r)$ and $f_k|_{\Pi^k}=0$, passing to
the limit as $s\to\infty$ we obtain
$$
(-1)^{i-k-1}\int\limits_{r_1}^{r_2}dr\int\limits_{B(a,r)}J(x,f)\,d{\cal H}^n=
\int\limits_{r_1}^{r_2}\int\limits_{\Sigma(a,R)}f_i\,df_{1}\wedge\ldots\wedge
\widehat{df_{i}} \wedge\ldots\wedge df_{n}\,.
$$
We divide both sides of this equation by $r_2-r_1$. Letting $r_2\to r_1$,
we see that (\ref{eqOGF}) holds for almost all $R>0$.

Thus, (\ref{eq1840}) is proved completely.

Fix $0<\alpha<1$ and choose $\omega_0$ such that
$$
\alpha\,\varepsilon_f(a,R)\le \int\limits_{\Sigma(a,R)}
(*d\omega)\,d{\cal H}^{n-1}
\left/\int\limits_{\Sigma(a,R)}|\omega-\omega_0|\,d{\cal H}^{n-1}\right.\,.
$$
By (\ref{eq1840}) we have
$$
\alpha\,\varepsilon_f(a,R)\le \left.\int\limits_{\Sigma(a,R)}J(x,f)\,
d{\cal H}^{n}\right/V(a,R)
$$
and setting $\alpha\to 1$ we obtain
$$
\varepsilon_f(a,R)\le \left.\int\limits_{\Sigma(a,R)}J(x,f)\,
d{\cal H}^{n}\right/V(a,R)\,.
$$

Since $|\nabla d(a,x)|\equiv 1$ we arrive at
$$
\int\limits_{\Sigma(a,R)}J(x,f)\,d{\cal H}^{n-1}=
\int\limits_{S(a,R)}J(x,f)\,d{\cal H}^{n-1}=V'(a,R)\ge
\varepsilon_f(a,R)\,V(a,R)\,,
$$
which is true for a.e. $R\in (0,\infty)$.

The function $t \mapsto V(a,t)$ is absolutely continuous on $[r,R]$ and
solving this differential inequality we have
\begin{equation}
\label{eqVR1R2}
{\rm exp}\left\{\int\limits_{r}^{R}\varepsilon_f(\tau)\,d\tau\right\}\,\le
{{V(a,R)}\over {V(a,r)}}\; .
\end{equation}

We put
$$
l(a,t)=\min_{x\in S(a,t)}|f(x)-f(a)|\, , \quad L(a,t)=
\max_{x\in S(a,t)}|f(x)-f(a)|\,.
$$

Now,
$$
V(a,r)\ge \Omega_n\,l(a,r)^n \,,\quad
V(a,R)\le \Omega_n\, L(a,R)^n \,,
$$
where $\Omega_n={\cal H}^n\left(B(0,1)\right)\,.$
Hence,
\begin{equation}
\label{eqlL}
{{V(a,R)}\over {V(a,r)}}\le {{L(a,R)^n}\over {l(a,r)^n}}\,.
\end{equation}
\medskip

We use the following bound for quasiconformal maps.
\smallskip

\begin{theo}{}
\label{theoDcond}
Let $f:{\bf R}^n\rightarrow {\bf R}^n$ be a $K$-quasiconformal map.
Then for every point $a\in  {\bf R}^n$ and an arbitrary $\rho,\,0<\rho<\infty$
the following estimate holds
\begin{equation}
\label{eqD}
\frac{L(a,\rho)}{l(a,\rho)}\leq D_*\,,
D_* =  D_*(K)= \exp(4K(K+1) \sqrt{K-1}) \, \, .
\end{equation}
\end{theo}
\medskip

For a proof see,  \cite[Theorem {\bf 14.8}]{AVV}.
\medskip

Using (\ref{eqD}), we see that
\be
\label{eq3.1'}
\frac{L(a,R)}{l(a,r)}= {{L(a,R)}\over{l(a,R)}}
\,{{l(a,R)}\over{L(a,r)}}\, {{L(a,r)}\over{l(a,r)}}
\leq D_*^2\,\frac{l(a,R)}{L(a,r)}\,,\quad D_*=D_*(K)\,.
\ee

We next show that if $\beta= K^{1/(n-1)}$, then for $0<r<R$
\begin{equation}
\label{eqRr}
{{l(a,R)}\over{L(a,r)}}\le {\, }\, \left( \frac{R}{r} \right)^{\beta}\,.
\end{equation}
Indeed,  because the case $l(a,R)\le L(a,r)$ is clear,
we may assume that $l(a,R)>L(a,r)$.
Applying the well-known formula for the $n$-capacity of
a spherical condenser in ${\bf R}^n$, we obtain
$$
\begin{array}{ll}
\omega_{n-1}\left(\ln \displaystyle{\frac{l(a,R)}{L(a,r)}}\right)^{1-n}
&={\rm cap}\,\left( \overline{B}(f(a),L(a,r)), B(f(a),l(a,R))\right)\\ \\
\quad &\ge {\rm cap}\, \left(f\overline{B}(a,r),fB(a,R)\right)\\ \\
\quad &\ge \frac{1}{K}\,{\rm cap}\,\left(\overline{B}(a,r),B(a,R)\right)\\ \\
\quad &={1\over K}\,\omega_{n-1}\left(\ln {R\over r}\right)^{1-n}\,,\\ \\
\end{array}
$$
where $\omega_{n-1}$ is the $(n-1)$-dimensional
area of the boundary of $B(0,1)$.

From this inequality we arrive at (\ref{eqRr}).
Thus, by (\ref{eq3.1'}) and (\ref{eqRr}) we obtain
\be
\label{eq3.2}
{{L(a,R)}\over{l(a,r)}} \le
{\, }\,D_*^2\,\left( \frac{R}{r} \right)^{\beta}\,.
\ee

Using (\ref{eq3.1'}), we have now
$$
{{V(a,R)}\over {V(a,r)}}\le D_*^{2n}\,
 \left( {{l(a,R)}\over {L(a,r)}} \right)^n\,,
$$
and further by (\ref{eqRr}),
\begin{equation}
\label{eq2232}
{{V(a,R)}\over {V(a,r)}}\le
{\, }\,D_*^{2n}\,\,\left( \frac{R}{r} \right)^{n\beta}\,.
\end{equation}

The estimates (\ref{eqVR1R2}) and (\ref{eq2232}) imply the following statement.
\medskip

\begin{theo}{}
\label{theokn}
Let $\Pi^k$ be a $K$-quasiplane in ${\bf R}^n$ with $1\le k\le n-2$.
Then for every point $a\in \Pi^k$ and arbitrary numbers $0<r<R<\infty$,
we have
\begin{equation}
\label{epsirR}
{\rm exp}\left\{\int\limits_{r}^{R}\varepsilon_f(a,\tau)\,d\tau\right\}\,
\le {\, }\,D_*^{2n}(K)\,\left( \frac{R}{r} \right)^{n\beta}\,, \quad
\beta= K^{1/(n-1)} \, ,
\end{equation}
for every $K$-quasiconformal mapping $f:{\bf R}^n\to{\bf R}^n$ such that
$\Pi^k=f^{-1}(\Pi^k_0)$.
\end{theo}
\bigskip

\bigskip

\cc
\section{Proof of Main Theorem}{}

Fix a $K$-quasiconformal map $f: {\bf R}^n \to {\bf R}^n$ with
$f(\Pi^k) = \Pi^k_0 \, .$
We shall estimate $\varepsilon_f(a,R)$ by $\eta(\Sigma(a,R))$.
Setting 
$$
\omega_0=\sum_{i=k+1}^n(-1)^{i-k-1}\alpha_i\,df_1\wedge\ldots\wedge\widehat{df_i}\wedge
\ldots\wedge df_n\,, \quad
\varphi =
\displaystyle \sum_{i=k+1}^n
\left|df_1\wedge\ldots\wedge\widehat{df_i}\wedge\ldots \wedge df_n\right|^2 \, ,
$$
where $\alpha_i$ are arbitrary constants, we have
\begin{equation}
\label{eqvarepKR}
\varepsilon_f(a,R)\ge \int\limits_{\Sigma(a,R)}*d\omega
\left/\int\limits_{\Sigma(a,R)}|\omega-\omega_0|\,d{\cal H}^{n-1}\right.\,.
\end{equation}

Observe that
$$
\begin{array}{ll}
&\displaystyle\int\limits_{\Sigma(a,R)}
|\omega-\omega_0|\,d{\cal H}^{n-1}\\ \\
&\le
\displaystyle\int\limits_{\Sigma(a,R)}\displaystyle\sqrt{\sum_{i=k+1}^n
(f_i-\alpha_i)^2}\displaystyle\sqrt{\displaystyle\sum_{i=k+1}^n
\left|df_1\wedge\ldots\wedge\widehat{df_i}\wedge\ldots \wedge df_n\right|^2}\,d{\cal H}^{n-1}\\ \\
&\le
\left(\displaystyle\int\limits_{\Sigma(a,R)}\left(\displaystyle\sum_{i=k+1}^n
(f_i-\alpha_i)^2\right)^{n/2}
d{\cal H}^{n-1}\right)^{1/n}
\left( \displaystyle\int\limits_{\Sigma(a,R)}
\left( \varphi \right)^{n/(2(n-1))}\,d{\cal H}^{n-1} \right)^{(n-1)/n}\,.\\ \\
\end{array}
$$

It is not difficult to check that
\begin{equation}
\label{eqphink}
\left( \varphi \right)^{n/(2(n-1))}\le
c_1(n,k)\,\|f'(x)\|^n\,,
\end{equation}
where $c_1(n,k)$ is a constant and, in fact,
\begin{equation}
\label{eqc3nK}
c_1(n,k)=(n-k)^{n/(2(n-1))}\,.
\end{equation}
Indeed, since the form $df_1\wedge\ldots\wedge\widehat{df_i}\wedge\ldots \wedge df_n$ is
simple,
$$
\begin{array}{ll}
&|df_1\wedge\ldots\wedge\widehat{df_i}\wedge\ldots \wedge df_n|^2\\ \\
&=\langle
df_1\wedge\ldots\wedge\widehat{df_i}\wedge\ldots \wedge df_n,\,
df_1\wedge\ldots\wedge\widehat{df_i}\wedge\ldots \wedge df_n\rangle\\ \\
&\le\prod_{s=1\atop{s\ne i}}^n|df_s|^2\,.\\ \\
\end{array}
$$
Using the inequality between geometric and arithmetic means we obtain
$$
\left(\prod_{s=1\atop{s\ne i}}^n|df_s|^2\right)^{1\over{n-1}}\le {1\over{n-1}}
\sum_{s=1\atop{s\ne i}}^n|df_s|^2
$$
and hence,
\begin{equation}
\label{eqdf2wed}
|df_1\wedge\ldots\wedge\widehat{df_i}\wedge\ldots \wedge df_n|^{2/(n-1)}\le
{1\over{n-1}}\sum_{s=1\atop{s\ne i}}^n
|df_s|^2\,.
\end{equation}
Further, for an arbitrary vector $h\in {\bf R}^n$,
$$
|f'(x)h|^2=\sum_{s=1}^n\langle \nabla f_s, h\rangle^2\,,
$$
and setting $h=\left.\nabla f_i\right/|\nabla f_i|$, we find
\begin{equation}
\label{eq1556}
|\nabla f_i|^2+ {1\over{|\nabla f_i|^2}}\sum_{s=1\atop s\ne i}^n\langle
\nabla f_s,\nabla f_i\rangle^2\le \|f'\|^2\,.
\end{equation}
Thus, from (\ref{eqdf2wed}) it follows
$$
|df_1\wedge\ldots\wedge\widehat{df_i}\wedge\ldots \wedge df_n|\le \|f'(x)\|^{n-1}
$$
and (\ref{eqphink}) follows. Now (\ref{eqphink}) yields
$$
\int\limits_{\Sigma(a,R)}\left( \varphi \right)^{n/(2(n-1))}
d{\cal H}^{n-1}\le
c_1(n,k)\int\limits_{\Sigma(a,R)}\|f'(x)\|^n\,d{\cal H}^{n-1}\,.
$$

Next we use the inequality
$$
\left({1\over N}\sum_{i=1}^n |a_i|^{t_1}\right)^{1/t_1}\le
\left({1\over N}\sum_{i=1}^n |a_i|^{t_2}\right)^{1/t_2}\quad (t_1\le t_2)
$$
(see, for example, \cite[\S 16 Chapter I]{BB}).

We have
$$
\left({1\over{n-k}}\sum_{i=k+1}^n |f_i-\alpha_i|^2\right)^{1/2}\le
\left({1\over{n-k}}\sum_{i=k+1}^n |f_i-\alpha_i|^n\right)^{1/n}
$$
and hence,
$$
\left(\sum_{i=k+1}^n |f_i-\alpha_i|^2\right)^{1/2}\le c_2(n,k)
\left(\sum_{i=k+1}^n |f_i-\alpha_i|^n\right)^{1/n}
$$
with $c_2(n,k)=(n-k)^{(n-2)/(2n)}$. Now we see from the definition
(\ref{f11}) of $\eta(\Sigma(a,R))$ that
$$
\begin{array}{ll}
\displaystyle\int\limits_{\Sigma(a,R)}\left(\displaystyle\sum_{i=k+1}^n
(f_i-\alpha_i)^2\right)^{n/2}d{\cal H}^{n-1}
&\le c^n_2(n,k)\displaystyle\sum_{i=k+1}^n \displaystyle\int\limits_{
\Sigma(a,R)}|f_i-\alpha_i|^n\,d{\cal H}^{n-1} \\ \\
\quad&\le c^n_2(n,k)\eta^{-n}(\Sigma(a,R))\displaystyle\sum_{i=k+1}^n
\displaystyle\int\limits_{\Sigma(a,R)}
\left|\nabla_S f_i\right|^n\,d{\cal H}^{n-1}\,. \\ \\
\end{array}
$$

 However, (\ref{eq1556}) implies that $|\nabla_S f_i|\le \|f'\|$ for
$i=1,\ldots,n$ and hence,
$$
\sum_{i=k+1}^n |\nabla_S f_i|^n \le (n-k)\,\|f'\|^n\,.
$$
Thus, we arrive at the estimate
$$
\int\limits_{\Sigma(a,R)}\left(\displaystyle\sum_{i=k+1}^n
(f_i-\alpha_i)^2\right)^{n/2}d{\cal H}^{n-1}
\le c_3(n,k)\,\eta^{-n}(\Sigma(a,R))\displaystyle\int\limits_{\Sigma(a,R)}
\|f'(x)\|^n\,d{\cal H}^{n-1}\,.
$$
where
$$
c_3(n,k)=(n-k)\,c_2(n,k)=(n-k)^{n/2}\,.
$$

Substituting these estimates into (\ref{eqvarepKR}), we have
the inequality
$$
\varepsilon_f(a,R)\ge (n-k)^{-1}\,\eta(\Sigma(a,R))\int\limits_{\Sigma(a,R)}J(x,f)\,
d{\cal H}^{n-1}
\left/\int\limits_{\Sigma(a,R)}\|f'(x)\|^n\,d{\cal H}^{n-1}\right.\,.
$$
Thus, we obtained
\begin{equation}
\label{eqvarela}
\varepsilon_f(a,R)\ge {{1}\over{K_O(f)}}\,\eta(\Sigma(a,R))\,
\end{equation}
and using Theorem \ref{theokn}, we prove Theorem \ref{theoknl}.
\medskip

{\bf Acknowledgements.} The authors are indebted to
Istvan Prause for a number of very useful remarks and 
to the referee for a set of detailed suggestions for
the improvement of the text.

\medskip

\bigskip


\bigskip
 \noindent
 {\bf Martio }:\\
 Department of Mathematics and Statistics\\
 University of Helsinki\\
 00014 Helsinki\\
 FINLAND\\
 Email: {\tt martio@cc.helsinki.fi}\\
 Fax: 358-9-19151400\\

 \medskip

 \noindent
 {\bf Miklyukov}:\\
Mathematics Department\\
Volgograd State University\\
2 Prodolnaya 30\\
Volgograd 400062 \\
RUSSIA\\
 Email: {\tt miklyuk@hotmail.com}\\
 Fax: 8442-471608
 \medskip

 \noindent
 {\bf Vuorinen}:\\
 Department of Mathematics\\
 University of Turku\\
 20014 Turku\\
 FINLAND\\
 Email: {\tt vuorinen@utu.fi}\\
 Fax: 358-2-3336595\\


\small

\begin{thebibliography}{AVV}

\bibitem{A} {\sc L.~Ahlfors,}
Quasiconformal reflections, Acta Math. {\bf 109} (1963), 291--301.

\bibitem{AVV}
{\sc G.D.~Anderson, M.K.~Vamanamurthy, and M.K.~Vuorinen,} Conformal
Invariants, Inequalities, and Quasiconformal Maps, Canadian Math. Soc., Ser.
of Monographs and Advanced texts, A Wiley - Interscience Publication, John
Wiley\&Sons, Inc., New-York -- Chichester -- Weinheim -- Brisbane -- Singapore
-- Toronto, 1997.

\bibitem{B}
{\sc C.~Bandle,}  Isoperimetric inequalities and applications, Pitman
Advanced Publishing, Boston-London-Melbourne, 1980.

\bibitem{BB}
{\sc E.F.~Beckenbach and  R.~Bellman}, Inequalities, Springer-Verlag,
Berlin -- G\"ottingen -- Heidelberg, 1961.


\bibitem{BH}
{\sc M. Bonk and J. Heinonen,} Smooth quasiregular mappings with branching,
  Publ. Math. Inst. Hautes £tudes Sci.  No. 100 (2004), 153--170.


\bibitem{EG}
{\sc L.C.~Evans and R.F.~Gariepy},
{Measure Theory and Fine Properties of Functions},
Studies in Advanced Mathematics, CRC PRESS,
Boca Raton -- New York -- London -- Tokio, 1992.

\bibitem{FMMVW}
{\sc D.~Franke, O.~Martio, V.M.~Miklyukov, M.~Vuorinen, and R.~Wisk,}
Quasiregular mappings and ${\cal WT}-$classes of differential forms
on Riemannian manifolds, Pacific J. Math. 202 (2002), 73--92.

\bibitem{G}
{\sc F.W.~Gehring,} Uniform domains and the ubiquitous quasidisk,
Jahresber. Deutsch. Math.-Verein, {\bf 89}  (1987), 88--103.

\bibitem{Gr}
{\sc M. Gromov, } Pseudoholomorphic curves in symplectic manifolds,
Invent. Math. 82 (1985),  307--347.

\bibitem{HKM}
{\sc J.~ Heinonen, T.~ Kilpel\"ainen and O.~ Martio,} Nonlinear potential
theory of degenerate elliptic equations, Clarendon Press, 1993.

\bibitem{Kr}
{\sc S.L.~Krushkal,} Quasiconformal mirrors,
Sibirsk. Mat. Zh. 40 (1999), 880--892.

\bibitem{MMPV}
{\sc O.~Martio,  V.~Miklyukov, S.~Ponnusamy, and M.~Vuorinen,}
On Some Properties of Quasiplanes, Results Math. 42 (2002),
107--113.

\bibitem{MV}
{\sc P. Mattila and M.Vuorinen,} Linear approximation property,
Minkowski dimension and quasiconformal spheres,
 J. London Math. Soc. (2) 42 (1990), 249--266.


\bibitem{Mik}
{\sc V.M.~Miklyukov,} On quasiconformally flat surfaces in Riemannian
manifolds, Izv. RAN (ser. math.), 67 (2003), 83--106.


\bibitem{Simons:1968}
\textsc{J.~Simons,}
{ Minimal varieties in Riemannian manifolds,}
{Ann. Math.}, 88 (1968), No.~1,
62--105.

\bibitem{VVW} \textsc{J. V\"ais\"al\"a, M. Vuorinen and H. Wallin,}
Thick sets and quasisymmetric maps, 
Nagoya Math. J. 135 (1994), 121-148.




\bibitem{Vu2}
{\sc  M.~Vuorinen,}
Geometric properties of quasiconformal maps and special functions. I.
Quasiconformal maps and spheres. (English. English, Polish summary)
Bull. Soc. Sci. Lett. \L\'od\'z S\'er. Rech. D\'eform. 24 (1997), 7--22.
(see also Errata: "Geometric properties of quasiconformal maps and
special functions. I, II, III [Bull. Soc. Sci. Lett. \L\'od\'z S\'er. Rech.
D\'eform. 24 (1997), 7--22; ibid., 23--35; ibid., 37--58; ].
Bull. Soc. Sci. Lett. \L\'od\'z S\'er. Rech. D\'eform.  26  (1998).



\bibitem{MMV3}
{\sc M.~Vuorinen, O.~Martio, and V.M.~Miklyukov,} On geometric structure
of quasiplanes, Proc. of dept. of mathematical analysis and
function theory, Volgograd State Univ., 2002, 21--31, ISBN 5-85534-531-9.


\end{thebibliography}
\end{document}